\documentclass[a4paper,12pt]{article}

\usepackage{amssymb}

\usepackage{amsmath}
\usepackage{amsfonts}
\usepackage{amscd}
\usepackage[english]{babel}

\title{The space of leftinvariant orthogonal almost complex structures on 6-dimensional Lie groups}
\author{Natalia Daurtseva\footnote{natali0112@ngs.ru}\\ Russia, Kemerovo State University\footnote{The work was supported by RFBR, 12-01-00873-a and by Russian President Grant supporting scientific schools SS-544.2012.1.}}
\date{ }
\begin{document}
\maketitle

\begin{abstract}
The space $\mathcal{Z}$ of leftinvariant orthogonal almost complex structures,
keeping the orientation, on 6-dimensional Lie groups is researched.
To get explicit view of this space elements the isomorphism of $\mathcal{Z}$ and $\mathbb{C}P^3$ is used.
The explicit formula for almost complex structure as composition of rotations is found.
\end{abstract}

\subsection*{1. Introduction}
Let $G$ -- 6-dimensional Lie group. Any leftinvariant almost complex structure
on Lie group is identified with Lie algebra $\mathfrak{g}$ endomorphism $I$, such that $I^2=-1$.
If a leftinvariant metric is fixed on the Lie group, then one can define subset of all leftinvariant
almost complex structures, which keep the metric. The space of this sort of almost complex structures, with additional property to keep orientation
is homogeneous space $\mathcal{Z}=SO(6)/U(3)$ and it is isomorphic to $\mathbb{C}P^3$.

In different problems with Hermitian structures on manifolds, and particulary on Lie groups, sometimes we need in explicit formulas
for almost complex structures, instead implicit condition $I^2=-1$. There exists only one orthogonal leftinvariant almost complex structure in dimension 2.
These structures form 2-parametric family
equal to $S^2$ in dimension 4. In 6-dimensional case, we have no explicit formulas of this sort of structures.

The main result of the article is explicit description of structures in $\mathcal{Z}$, via homogeneous coordinates in $\mathbb{C}P^3$,
and via six angles, as composition of rotations.

\subsection*{2. Notations and definitions}

Let $(M^{2n},g)$ is Riemannian manifold of class $C^{\infty}$.

\vspace{1.5mm}
\textbf{Definition 1.} {\it Almost complex structure} on $M$ is smooth endomorphism field $J_x:T_xM\longrightarrow T_xM$,
such that $J_x^2=-Id_x$, where $Id_x$ - identical endomorphism of $T_xM$, $\forall x\in M$.
If almost complex structure keep $g$:
$$
g(JX,JY)=g(X,Y),\qquad\forall X,Y\in\Xi(M),
$$
then $J$ is called {\it orthogonal} with respect to metric $g$, and pair $(g,J)$ on $M^{2n}$ is called {\it almost Hermitian structure.}
Manifold with almost Hermitian structure is called {\it almost Hermitian.}
\vspace{1.5mm}

\vspace{1.5mm}
\textbf{Definition 2.} Almost complex structure $J$ on $M$ is called {\it associated with skew-symmetric 2-form} $\omega$, if
$$
\omega(JX,JY)=\omega(X,Y),\qquad\forall X,Y\in\Xi(M)
$$
If $(g,J)$ is almost Hermitian structure on $M$, then such 2-form $\omega$ is defined by formula
$$
\omega(X,Y)=g(JX,Y)
$$
and is called {\it fundamental 2-form.}
 \vspace{1.5mm}

\textsc{Remark 1.} We will identify almost complex structure on Riemannian manifold with its fundamental form,  when it will  be necessary.
\vspace{1.5mm}

\vspace{1.5mm}
\textsc{Remark 2.} Every almost complex structure on the manifold defines orientation on it [3].
\vspace{1.5mm}

Let $M=G$ is 6-dimensional Lie group, $g$ -- some leftinvariant metric on $G$. Leftinvariant almost complex structure $I$
on $G$ is uniquely defined by its restriction on $\mathfrak{g}$, by linear endomorphism:
\begin{equation}
I:\mathfrak{g}\longrightarrow\mathfrak{g},\qquad I^2=-1
\end{equation}
On another hand, linear endomorphism (1) defines leftinvariant almost complex structure on Lie group $G$.
In case of leftinvariant almost complex structures on Lie groups we will use the same notations to these structures and there restrictions on Lie algebras.
So, the set of all $g$-orthogonal leftinvariant almost complex structures on Lie group $G$ is the set of all linear endomorphisms:
$$
\begin{array}{c}
I:\mathfrak{g}\longrightarrow\mathfrak{g},\qquad I^2=-1,\\
g(IX,IY)=g(X,Y),\qquad\forall X,Y\in\mathfrak{g}
\end{array}
$$
Set $\mathcal{Z}$ of all these structures, keeping orientation, is homogeneous space $SO(6)/U(3)$. Really [3],
fix some orthonormal basis $(e)=(e_1,\dots,e_6)$ of Lie algebra $\mathfrak{g}$. Let $I_0$ is almost complex structure $I_0e_1=-e_4$, $I_0e_2=-e_5$, $I_0e_3=-e_6$,
then any another orthogonal leftinvariant almost complex structure $I$ is equal to $SI_0S^{-1}$, for some $S\in SO(6)$,
i.e. group $SO(6)$ acts transitively on $\mathcal{Z}$. Isotropy subgroup of $I_0$ consists of orthogonal matrixes, which commute
with $I_0$, i.e. coincide with $U(3)$. Thus $\mathcal{Z}=SO(6)/U(3)$.

\subsection*{3. Diffeomorphism of $\mathcal{Z}$ in $\mathbb{C}P^3$}
It is known [4], that $SO(6)/U(3)$ is diffeomorphic to $\mathbb{C}P^3$. Construct this diffeomorphism explicitly.
Let $(v^0,v^1,v^2,v^3)$ is unitary basis $V=\mathbb{C}^4$. Then $\Lambda^2\mathbb{C}^4$ is identified with  $\mathfrak{g}^*$ by:
\begin{equation}
\begin{array}{ll}
2v^0\wedge v^1=e^1+ie^4, & 2v^2\wedge v^3=e^1-ie^4,\\
2v^0\wedge v^2=e^2+ie^5, & 2v^3\wedge v^1=e^2-ie^5,\\
2v^0\wedge v^3=e^3+ie^6, & 2v^1\wedge v^2=e^3-ie^6,\\
\end{array}
\end{equation}
here $e^i$ is co-vector, dual to $e_i$. Almost complex structure $I$ on $\mathfrak{g}$ defines almost complex structure on the $\mathfrak{g}^*$,
by the standard way: $(Ie^i)(X)=e^i(IX)$. We will use the same notation to this structure.
Every almost complex structure defines decomposition of complexification
$\mathfrak{g}^{*\mathbb{C}}$ into the direct sum of proper subspaces:
$$
\begin{array}{l}
\mathfrak{g}^{*1,0}=\{\alpha\in\mathfrak{g}^{*\mathbb{C}}:I\alpha=i\alpha\}=\{\varphi-iI\varphi:\varphi\in\mathfrak{g}^*\}\\
\mathfrak{g}^{*0,1}=\{\alpha\in\mathfrak{g}^{*\mathbb{C}}:I\alpha=-i\alpha\}=\{\varphi+iI\varphi:\varphi\in\mathfrak{g}^*\}
\end{array}
$$
Vice versa, every almost complex structure is uniquely defined by $\mathfrak{g}^{*1,0}$.

For standard structure $I=I_0$, proper subspace, corresponding to $i$
is $\mathfrak{g}^{*1,0}=\mbox{span}_{\mathbb{C}}\{e^1+ie^4,e^2+ie^5,e^3+ie^6\}$, thus by (2) space $\mathfrak{g}^{*1,0}$
is identified with $V_{v^0}=\{v^0\wedge v,\ v\in V\}\subset\Lambda^2V$.

Arbitrary almost complex structure $I$ has the same property as $I_0$ in another basis $(e')=(e)S$, where $S\in SO(6)$.
As $SU(4)\longrightarrow SO(6)$ is double covering, then for arbitrary almost complex structure the subspace in $\Lambda^2V$, corresponding to
$\mathfrak{g}^{*1,0}$ is:
$$
V_u=\{u\wedge v:v\in V\}\subset\Lambda^2V
$$
So,  $\mathcal{Z}$ is in one to one correspondence with point with $\mathbb{C}P^3$ by the following way:
$$
I\in\mathcal{Z}\longrightarrow\mathfrak{g}^{*1,0}\longrightarrow V_u\longrightarrow [u]\in\mathbb{C}P^3
$$
Let define the following almost complex structures:
$$
\begin{array}{c}
I_0e^1=-e^4, I_0e^2=-e^5, I_0e^3=-e^6;\\
I_1e^1=-e^4, I_1e^2=e^5, I_1e^3=e^6;\\
I_2e^1=e^4, I_2e^2=-e^5, I_2e^3=e^6;\\
I_3e^1=e^4, I_3e^2=e^5, I_3e^3=-e^6;\\
\end{array}
$$
All these ones give the same orientation, keep metric $g$, and correspond to the points
$[1,0,0,0], [0,1,0,0], [0,0,1,0], [0,0,0,1]\in\mathbb{C}P^3$.
Fundamental 2-forms, corresponding to these structures are:
$$
\begin{array}{c}
\omega_0=e^1\wedge e^4+e^2\wedge e^5+e^3\wedge e^6;\\
\omega_1=e^1\wedge e^4-e^2\wedge e^5-e^3\wedge e^6;\\
\omega_2=-e^1\wedge e^4+e^2\wedge e^5-e^3\wedge e^6;\\
\omega_3=-e^1\wedge e^4-e^2\wedge e^5+e^3\wedge e^6.
\end{array}
$$
Let $z=[z^0,z^1,z^2,z^3]$ and $u=[u^0,u^1,u^2,u^3]$ are arbitrary points into $\mathbb{C}P^3$.
One can join these points by 'edge', consisting of points
$\alpha z+\beta u$, where $\alpha,\beta\in\mathbb{C},\quad |\alpha|+|\beta|\neq 0$, let's denote it as $\mathcal{E}_{zu}$.
It is not difficult to show, that 'edge' $\mathcal{E}_{zu}=\{[\alpha,\beta]:\alpha,\beta\in\mathbb{C},\quad |\alpha|+|\beta|\neq 0\}=\mathbb{C}P^1=S^2$.
By analogy three arbitrary points, which are not on the same 'edge' define the face $\mathcal{F}=\mathbb{C}P^2$.

Thus, projective space $\mathbb{C}P^3$ one may visualize
as solid tetrahedron with vertexes
$\omega_0=[1,0,0,0]$, $\omega_1=[0,1,0,0]$, $\omega_2=[0,0,1,0]$, $\omega_3=[0,0,0,1]$, 'edges' $\mathcal{E}_{ij}\cong\mathbb{C}P^1$
and 'faces' $\mathcal{F}_i\cong\mathbb{C}P^2$ (edge $\mathcal{E}_{ij}$ joins forms $\omega_i$ and $\omega_j$,
face $\mathcal{F}_i$ is 'opposite' to vertex $\omega_i$).

\vspace{1.5mm}
\textbf{Lemma 1.} {\it Fundamental 2-form $\omega\in\mathcal{E}_{01}$ is of view:}
\begin{equation}
\omega=e^1\wedge e^4+r(e^2\wedge e^5+e^3\wedge e^6)
+u(e^2\wedge e^3+e^6\wedge e^5)+x(e^2\wedge e^6+e^5\wedge e^3)
\end{equation}
{\it where $r^2+u^2+x^2=1$.}
\vspace{1.5mm}

{\bf Proof:} For arbitrary form $\omega\in E_{01}$ there exists numbers
$s,c_1,c_2\in\mathbb{R}$, $s^2+c_1^2+c_2^2=1$, such that corresponding almost complex structure is:
$I=sI_0+(c_1+ic_2)I_1$:
$$
s[1,0,0,0]+(c_1+ic_2)[0,1,0,0]=[s, c_1+ic_2,0 ,0].
$$
Corresponding space
$$V_{[s, c_1+ic_2,0 ,0]}=\{sv^0+(c_1+ic_2)v^1\wedge u, u\in V\}.$$

As $(sv^0+(c_1+ic_2)v^1)\wedge v^0=(c_1+ic_2)v^1\wedge v^0=-\frac12(c_1+ic_2)(e^1+ie^4)=\frac12(-c_1e^1+c_2e^4+i(-c_2e^1-c_1e^4))$,
then for almost complex structure $I\in\mathcal{E}_{01}$:
$$
I(-c_1e^1+c_2e^4)=c_2e^1+c_1e^4.
$$
By analogy, equalities
$$
(sv^0+(c_1+ic_2)v^1)\wedge v^1=\frac12(se^1+ise^4);
$$
$$
(sv^0+(c_1+ic_2)v^1)\wedge v^2=sv^0\wedge v^2+(c_1+ic_2)v^1\wedge v^2=
$$
$$
=\frac12(s(e^2+ie^5)+(c_1+ic_2)(e^3-ie^6))=
$$
$$
=\frac12(se^2+c_1e^3+c_2e^6+i(se^5+c_2e^3-c_1e^6));
$$
$$
(sv^0+(c_1+ic_2)v^1)\wedge v^3=sv^0\wedge v^3+(c_1+ic_2)v^1\wedge v^3=
$$
$$
\frac12(s(e^3+ie^6)+(c_1+ic_2)(-e^2+ie^5))=
$$
$$
=\frac12(se^3-c_1e^2-c_2e^5+i(se^6-c_2e^2+c_1e^5))
$$
give:
$$
I(se^1)=-se^4
$$
$$
I(se^2+c_1e^3+c_2e^6)=-(se^5+c_2e^3-c_1e^6)
$$
$$
I(se^3-c_1e^2-c_2e^5)=-(se^6-c_2e^2+c_1e^5)
$$
Then $\omega\in\mathcal{E}_{01}$, corresponding to $I$ is of view:
$$
\omega=(-c_1e^1+c^2e^4)\wedge (-c_2e^1-c^1e^4)+se^1\wedge se^4+(se^2+c_1e^3+c_2e^6)\wedge (se^5+c_2e^3-c_1e^6)+
$$
$$
+(se^3-c_1e^2-c_2e^5)\wedge (se^6-c_2e^2+c_1e^5)=
$$
$$
=e^1\wedge e^4+(s^2-c_1^2-c_2^2)e^2\wedge e^5+2sc_2e^2\wedge e^3-
2sc_1e^2\wedge e^6-2sc_1e^5\wedge e^3-
$$
$$
-2sc_2e^5\wedge e^6+(s^2-c_1^2-c_2^2)e^3\wedge e^6
$$
Suppose $r=s^2-c_1^2-c_2^2=2s^2-1$, $u=2sc_2$, $x=-2sc_1$, then $\omega$ is given by (3).
\begin{flushright}$\Box$\end{flushright}

\textbf{Theorem 1.}
{\it
Let $[1,a,b,c]\in\mathbb{C}P^3$ is arbitrary point out of face $\mathcal{F}_0$. Leftinvariant almost complex structure $I\in\mathcal{Z}$,
corresponding to this point is in the neighborhood $U(I_0)=\{I\in\mathcal{Z}:1-II_0 - \mbox{ is inveresed}\}$ of structure $I_0$ and in the basis $(e)$
is given by matrix:}
\footnotesize
\begin{equation}
\setlength{\arraycolsep}{1pt}
\frac{1}{x}
\left(
\begin{array}{cccccc}
0 & 2\Im(\overline{a}b+c) & 2\Im(\overline{a}c-b) & x-2|b|^2-2|c|^2 & 2\Re(\overline{a}b-c) & 2\Re(\overline{a}c+b)\\
-2\Im(\overline{a}b+c) & 0 & 2\Im(\overline{b}c+a) & 2\Re(\overline{a}b+c) & x-2|a|^2-2|c|^2 & 2\Re(\overline{b}c-a)\\
-2\Im(\overline{a}c-b) & -2\Im(\overline{b}c+a) & 0 & 2\Re(\overline{a}c-b) & 2\Re(\overline{b}c+a) & x-2|a|^2-2|b|^2\\
2|b|^2+2|c|^2-x & -2\Re(\overline{a}b+c) & -2\Re(\overline{a}c-b) & 0 & 2\Im(\overline{a}b-c) & 2\Im(\overline{a}c+b)\\
-2\Re(\overline{a}b-c) & 2|a|^2+2|c|^2-x & -2\Re(\overline{b}c+a) & -2\Im(\overline{a}b-c) & 0 & 2\Im(\overline{b}c-a)\\
-2\Re(\overline{a}c+b) & -2\Re(\overline{b}c-a) & 2|a|^2+2|b|^2-x & -2\Im(\overline{a}c+b) & -2\Im(\overline{b}c-a) & 0
\end{array}
\right)
\end{equation}
\normalsize
{\it Let $I$ is arbitrary almost complex structure in the neighborhood $U(I_0)$, and $(I_{ij})_{i,j=1,..,6}$
is matrix of this structure in $(e)$.
Then coordinates of corresponding point $[1,a,b,c]$ of $\mathbb{C}P^3$ are following:}
\[
\begin{array}{c}
a=\frac{1}{1+I_{14}+I_{25}+I_{36}}((I_{35}-I_{26})+i(I_{23}-I_{56}))\\
b=\frac{1}{1+I_{14}+I_{25}+I_{36}}((I_{16}-I_{34})+i(I_{46}-I_{13}))\\
c=\frac{1}{1+I_{14}+I_{25}+I_{36}}((I_{24}-I_{15})+i(I_{12}-I_{45}))
\end{array} \eqno(5)
\]
{\it where $x=1+|a|^2+|b|^2+|c|^2$.
}
 \vspace{1.5mm}

{\bf Proof:} Let $a,b,c\in\mathbb{C}$ and $[1,a,b,c]$ is some point from $\mathbb{C}P^3$,
out of face $\mathcal{F}_0$. Then direct calculations, which analogous to the ones in lemma 1 proof gives corresponding form
$\omega$ and almost complex structure (4).

Let $I\in U(I_0)=\{I\in\mathcal{Z}:1-II_0  \mbox{ is inveresed }\}$. For this structure [2] the corresponding skew-symmetric operator $K$, which anticommutate with
$I_0$ is defined by:
$$
\begin{array}{c}
I=(1-K)I_0(1-K)^{-1}\\
K=(1-II_0)^{-1}(1+II_0)
\end{array}
$$
Let $K=\left(\begin{array}{cc}
A & B\\
B & -A
\end{array}
\right),$
where $A+iB=\left(\begin{array}{ccc}
0 & -c & b\\
c & 0 & -a\\
-b & a & 0
\end{array}
\right)$, where $I$ is defined by (4).

Vice versa, let $I=(I_{ij})$ ($I_{ij}=-I_{ji}$) is arbitrary almost complex structure from $U(I_0)$,
then formulas (5) are deduced directly from (4).
\begin{flushright}$\Box$\end{flushright}

\vspace{1.5mm}
\textbf{Corollary 1.} {\it
Let $I=[1,a,b,c]\in\mathbb{C}P^3$, then operator $K=(1-II_0)^{-1}(1+II_0)$ in $(e)$ is equal to:
$\left(\begin{array}{cc}
A & B\\
B & -A
\end{array}
\right),$
where $A+iB=\left(\begin{array}{ccc}
0 & -c & b\\
c & 0 & -a\\
-b & a & 0
\end{array}
\right)$
}
\vspace{1.5mm}

\textsc{Remark 3.} Formula (4) provide explicit description of orthogonal matrixes satisfied to condition $I^2=-1$.
 \vspace{1.5mm}

It is known [1], that $\mathbb{C}P^3/T^3$ is 3-dimensional solid tetrahedron. Let's define action of $T^3$ on
$\mathcal{Z}$.

If forms $p_+$ and $p_-$
get all possible values on edges $\mathcal{E}_{03}$ and $\mathcal{E}_{12}$     accordingly, then:
edges $\mathcal{E}_{p_+p_-}$ fill the whole 6-dimensional tetrahedron $\mathbb{C}P^3$.
$$
{\mathcal Z}=\bigcup_{
\begin{array}{c}
p_+\in \mathcal{E}_{03}\\
p_-\in \mathcal{E}_{12}
\end{array}}\mathcal{E}_{p_+p_-}.
$$

Fix 2-dimensional platform $\mathbb{D}_1=\langle e^1,e^4\rangle$, $\mathbb{D}_2=\langle e^2,e^5\rangle$, $\mathbb{D}_3=\langle e^3,e^6\rangle$.
Torus $T^3$
acts in $\mathbb{R}^6$ by rotations into these platforms.

We will think that two almost complex structures $I$ and $J$
will be in the same equivalence class in $\mathcal{Z}/T^3$, if $I=OJO^{-1}$, where
$O\in SO(6)$ is rotation into these platforms matrix. By another words, we will identify almost complex structures,
which differ each others on the rotations into these platforms only.
As $(\cos\alpha e+\sin\alpha f)\wedge(-\sin\alpha e+\cos\alpha f)=e\wedge f$, then tetrahedron vertexes
are invariant under this action.

Let use explicit view of edge $\mathcal{E}_{01}$ from lemma 1. Introduce spherical coordinates:
$$
\begin{array}{lc}
r=\sin\psi; &\\
u=\cos\varphi\cos\psi; & -\pi/2\leq\psi\leq\pi/2\\
v=\sin\varphi\cos\psi; & 0\leq\varphi\leq2\pi
\end{array}
$$
Then form $\omega\in\mathcal{E}_{01}$ is:
$$
\omega=e^1\wedge e^4+\sin\psi(e^2\wedge e^5+e^3\wedge e^6)+\cos\varphi\cos\psi(e^2\wedge e^3+e^6\wedge e^5)+
\sin\varphi\cos\psi(e^2\wedge e^6+e^5\wedge e^3)=
$$
$$
=e^1\wedge e^4+\sin\psi(e^2\wedge e^5+e^3\wedge e^6)+
\cos\psi(\cos\varphi e^2\wedge e^3+\cos\varphi e^6\wedge e^5+\sin\varphi e^2\wedge e^6+\sin\varphi e^5\wedge e^3).
$$
Let
$$
f^2=\cos\varphi e^2+\sin\varphi e^5; f^5=-\sin\varphi e^2+\cos\varphi e^5,
$$
then $f^2\wedge e^3+e^6\wedge f^5=\cos\varphi f^2\wedge e^3+\cos\varphi e^6\wedge f^5+\sin\varphi f^2\wedge e^6+\sin\varphi f^5\wedge e^3$,
and following the form $\omega\in\mathcal{E}_{01}/SO(2)$ is  $e^1\wedge e^4+\sin\psi(e^2\wedge e^5+e^3\wedge e^6)+\cos\psi(e^2\wedge e^3+e^6\wedge e^5)$.
In this case two almost complex structures on sphere-edge $\mathcal{E}_{01}$ are equivalent if they are on the same parallel of this sphere.

What happens in general case for the edge $\mathcal{E}_{p_+p_-}$? There exists [4] the following description of the form
$\omega\in\mathcal{Z}$. Let $J\in\mathcal{Z}$ is almost complex structure, corresponding to $\omega$ and
$\mathbb{D}=\langle e^1,e^2,e^4,e^5\rangle$. Fix $e^3$,
then $Je^3$ is orthogonal to $e^3$, and so $Je^3=ae^6+bf^1$, where $f^1\in\mathbb{D}$, $\|f^1\|=1$ and $a^2+b^2=1$.
Unit  1-form
$af^1-be^6$ is orthogonal to $e^3$ and $Je^3$, thus 2-form $\omega\in\mathcal{Z}$ is:
$$
\omega(P,a,b)=e^3\wedge(ae^6+bf^1)-f^4\wedge(af^1-be^6)+f^2\wedge f^5
$$
here $(f^1,f^2,f^4,f^5)=(e^1,e^2,e^4,e^5)S$, for suitable $S\in SO(4)$. It is known, that for $\Lambda^2\mathbb{D}$ one have decomposition:
$
\Lambda^2\mathbb{D}=\Lambda^2_+\mathbb{D}\oplus\Lambda^2_-\mathbb{D},
$
where
$$
\begin{array}{cc}
\Lambda^2_+\mathbb{D}=\{e^{14}+e^{25},e^{12}+e^{54},e^{15}+e^{42}\}\\
\Lambda^2_-\mathbb{D}=\{e^{14}-e^{25},e^{12}-e^{54},e^{15}-e^{42}\}
\end{array}
$$
are proper subspaces of operator * ($e^{ij}$ denotes 2-form $e^i\wedge e^j$).
This decomposition provide double covering
$SO(4)\longrightarrow SO(3)\times SO(3)$.
Thus arbitrary matrix $P\in SO(4)$ is represented by matrix
$\left(\begin{array}{cc}
P_+ & 0\\
0 & P_-
\end{array}\right)\in SO(6)$, where $P_+,P_-\in SO(3)$. For $\omega$ matrixes $P_+$ and $P_-$ have the clear geometric sense, namely
$p_+=\omega(P,1,0)\in\mathcal{E}_{03}$, $p_-=\omega(P,-1,0)\in\mathcal{E}_{12}$. When $P_+$ is changing, then form $p_+$ is moving along the edge
$\mathcal{E}_{03}$, and if $P_-$ is changing then form $p_-$ is moving along the edge $\mathcal{E}_{12}$. Therefore, accurate within
to rotations into the platforms
$\mathbb{D}_1, \mathbb{D}_2, \mathbb{D}_3$:
$$
\mathcal{E}_{03}/S^1=\{\sin\psi(e^{14}+e^{25})+\cos\psi(e^{12}+e^{54})+e^{36}:-\pi/2\leq\psi\leq\pi/2\}
$$
$$
\mathcal{E}_{12}/S^1=\{\sin\theta(e^{14}-e^{25})+\cos\theta(e^{12}-e^{54})-e^{36}:-\pi/2\leq\theta\leq\pi/2\}
$$
One can check that matrix
$\left(
\begin{array}{cccc}
\sin\frac{\psi+\theta}{2} & 0 & 0 & -\cos\frac{\psi+\theta}{2}\\
0 & \cos\frac{\psi-\theta}{2} & -\sin\frac{\psi-\theta}{2} & 0\\
0 & \sin\frac{\psi-\theta}{2} & \cos\frac{\psi-\theta}{2} & 0\\
\cos\frac{\psi+\theta}{2} & 0 & 0 & \sin\frac{\psi+\theta}{2}
\end{array}\right)$ from $SO(4)$
corresponds to the pair of matrices $(P_+,P_-)$, where
$P_+=
\left(
\begin{array}{ccc}
\sin\psi & 0 & 0\\
0 & \cos\psi & 0\\
0 & 0 & 1
\end{array}
\right)$,

$P_-=\left(
\begin{array}{ccc}
\sin\theta & 0 & 0\\
0 & \cos\theta & 0\\
0 & 0 & 1
\end{array}
\right).$
Suppose $a=\sin\varphi$, $b=\cos\varphi$, $\varphi\in [0,2\pi]$, thus we get:

\vspace{1.5mm}
\textbf{Theorem 2.} {\it Form $\omega\in\mathcal{Z}/\mathbb{T}^3$ is:
\begin{center}
$
\omega=e^3\wedge(\sin\varphi e^6+\cos\varphi(\sin\frac{\psi+\theta}{2}e^1+\cos\frac{\psi+\theta}{2}e^5))+
(\sin\varphi(\sin\frac{\psi+\theta}{2}e^1+\cos\frac{\psi+\theta}{2}e^5)-\cos\varphi e^6)\wedge (-\sin\frac{\psi-\theta}{2}e^2+\cos\frac{\psi-\theta}{2}e^4)+
(\cos\frac{\psi-\theta}{2}e^2+\sin\frac{\psi-\theta}{2}e^4)\wedge (-\cos\frac{\psi+\theta}{2} e^1+\sin\frac{\psi+\theta}{2}e^5)
$
\end{center}
where $\varphi,\psi,\theta\in [-\pi/2,\pi/2]$.}
 \vspace{1.5mm}

\vspace{1.5mm}
\textbf{Corollary 2.} {\it
Components $(J_{ij})_{i,j=1}^6$ of arbitrary almost complex structures $J\in\mathcal{Z}$ in $(e)$ are:}
\footnotesize
$$
\begin{array}{l}
J_{12}=\sin\frac{\psi+\theta}{2}\sin\frac{\psi-\theta}{2}(\cos\varphi_1\cos\varphi_2-\sin\varphi\sin\varphi_1\sin\varphi_2)+\cos\frac{\psi+\theta}{2}\cos\frac{\psi-\theta}{2}(\sin\varphi_1\sin\varphi_2-\sin\varphi\cos\varphi_1\cos\varphi_2)\\
J_{13}=\cos\varphi(\cos\frac{\psi-\theta}{2}\cos\varphi_1\cos\varphi_3-\sin\frac{\psi+\theta}{2}\sin\varphi_1\sin\varphi_3)\\
J_{14}=\sin\varphi\sin\frac{\psi+\theta}{2}\cos\frac{\psi-\theta}{2}+\cos\frac{\psi+\theta}{2}\sin\frac{\psi-\theta}{2}\\
J_{15}=\cos\frac{\psi+\theta}{2}\cos\frac{\psi-\theta}{2}(\sin\varphi\cos\varphi_1\sin\varphi_2+\sin\varphi_1\cos\varphi_2)-\sin\frac{\psi+\theta}{2}\sin\frac{\psi-\theta}{2}(\sin\varphi\sin\varphi_1\cos\varphi_2+\cos\varphi_1\sin\varphi_2)\\
J_{16}=-\cos\varphi(\sin\frac{\psi+\theta}{2}\sin\varphi_1\cos\varphi_3+\cos\frac{\psi-\theta}{2}\cos\varphi_1\sin\varphi_3)\\
J_{23}=\cos\varphi(\cos\frac{\psi+\theta}{2}\cos\varphi_2\sin\varphi_3+\sin\frac{\psi-\theta}{2}\sin\varphi_2\cos\varphi_3)\\
J_{24}=\sin\frac{\psi+\theta}{2}\sin\frac{\psi-\theta}{2}(\sin\varphi_1\cos\varphi_2+\sin\varphi\cos\varphi_1\sin\varphi_2)-\cos\frac{\psi+\theta}{2}\cos\frac{\psi-\theta}{2}(\sin\varphi\sin\varphi_1\cos\varphi_2+\cos\varphi_1\sin\varphi_2)\\
J_{25}=\sin\varphi\cos\frac{\psi+\theta}{2}\sin\frac{\psi-\theta}{2}+\cos\frac{\psi-\theta}{2}\sin\frac{\psi+\theta}{2}\\
J_{26}=\cos\varphi(\cos\frac{\psi+\theta}{2}\cos\varphi_2\cos\varphi_3-\sin\frac{\psi-\theta}{2}\sin\varphi_2\sin\varphi_3)\\
J_{34}=\cos\varphi(\sin\frac{\psi+\theta}{2}\cos\varphi_1\sin\varphi_3+\cos\frac{\psi-\theta}{2}\sin\varphi_1\cos\varphi_3)\\
J_{35}=\cos\varphi(\cos\frac{\psi+\theta}{2}\sin\varphi_2\sin\varphi_3-\sin\frac{\psi-\theta}{2}cos\varphi_2\cos\varphi_3)\\
J_{36}=\sin\varphi\\
J_{45}=\cos\frac{\psi+\theta}{2}\cos\frac{\psi-\theta}{2}(\cos\varphi_1\cos\varphi_2-\sin\varphi\sin\varphi_1\sin\varphi_2)-
\sin\frac{\psi+\theta}{2}\sin\frac{\psi-\theta}{2}(\sin\varphi\cos\varphi_1\cos\varphi_2-\sin\varphi_1\sin\varphi_2)\\
J_{46}=\cos\varphi(\cos\frac{\psi-\theta}{2}\sin\varphi_1\sin\varphi_3-\sin\frac{\psi+\theta}{2}\cos\varphi_1\cos\varphi_3)\\
J_{56}=-\cos\varphi(\cos\frac{\psi+\theta}{2}\sin\varphi_2\cos\varphi_3+\sin\frac{\psi-\theta}{2}\cos\varphi_2\sin\varphi_3),\ \mbox{\textit{ãäå}}\ \varphi,\psi,\theta\in [-\pi/2,\pi/2],\varphi_1,\varphi_2,\varphi_3\in [0,2\pi]
\end{array}
$$
\normalsize
\vspace{1.5mm}
\textsc{Remark 4.} Corollary 2 gives explicit description of arbitrary, leftinvariant almost complex structure
on 6-dimensional Lie group. As every orthogonal almost complex structure in some orthonormal basis acts as standard, then we can only
show how to get this suitable orthonormal basis.

Denote:

$R_{\varphi_i}$ is rotation into the platform $\mathbb{D}_i$ by angle $\varphi_i$,

$R_{\alpha}$ is rotation by angle $\alpha$
($\alpha=\frac{\psi+\theta}{2}$) into the platform
$\langle e_1,e_5\rangle$,

$R_{\beta}$ is the rotation by the angle $\beta$ ($\beta=\frac{\psi-\theta}{2}$) into the platform $\langle e_2,e_4\rangle$,

$R_{\varphi}$ is rotation by angle $\varphi$ into the platform $\langle e_1,e_6\rangle$.
Then we can get the suitable basis by the following composition of rotations:
$$
R_{\varphi}\circ R_{\alpha}\circ R_{\beta}\circ R_{\varphi_3}\circ R_{\varphi_2}\circ R_{\varphi_1}.
$$
\subsection*{4. Example} Let consider the group $G=SU(2)\times SU(2)$. One can accentuate
two special classes of almost complex structures on this Lie group.
Firstly it is  Calabi-Echkmann structures (integrable structures on the product of odd-dimensional spheres), and, secondly,
orthogonal leftinvariant structures,
which give maximum of Nijenhuis tensor norm ("maximum non-integrable" almost complex structures).

Calabi-Eckmann structures arise from Hopf fibration
$S^3\longrightarrow^{S^1}\mathbb{C}P^1$. Really on the product of 3-spheres it provides fibration with
complex base and complex fiber
$$S^3\times S^3\longrightarrow^{S^1\times S^1}\mathbb{C}P^1\times\mathbb{C}P^1.$$
Existence  of holomorphic transition functions allows to construct complex structure on the total space $S^3\times S^3$.
It is known [6,9], that the complex structures describe the whole class of leftinvariant complex structures on Lie group $SU(2)\times SU(2)$.
Fix standard basis  $(e_1,e_2,e_3,e_4,e_5,e_6)$ of Lie algebra
$\mathfrak{su}(2)\times\mathfrak{su}(2)=\mathbb{R}^3\times\mathbb{R}^3$. In this basis integrable structure, which give the same orientation
as $I_0, I_1, I_2, I_3$ is of view $Ie_1=-e_4$, $Ie_2=-e_3$, $Ie_5=e_6$.
Obviously, all other leftinvariant orthogonal complex structures are in the orbit of $I$ under the action of $SO(3)\times SO(3)$.
Stabilizer of this action is Lie group $SO(2)\times SO(2)$, so the set of integrable almost complex structures
is 4-dimensional subset $S^2\times S^2=SO(3)\times SO(3)/SO(2)\times SO(2)$.
By technical difficulties we will show not all this set in our model, but only some points of it.
Consider the following forms, corresponding to some integrable almost complex structures:
\begin{center}
$
\begin{array}{ll}
e^{14}\pm e^{23}\mp e^{56}\in\mathcal{E}_{01} &
-e^{14}\pm e^{23}\pm e^{56}\in\mathcal{E}_{23}\\
\pm e^{12}\mp e^{45}+e^{36}\in\mathcal{E}_{03} &
\pm e^{12}\pm e^{45}-e^{36}\in\mathcal{E}_{12}\\
e^{25}\pm e^{46}\mp e^{13}\in\mathcal{E}_{02} &
-e^{25}\pm e^{46}\pm e^{13}\in\mathcal{E}_{13}
\end{array}
$
\end{center}
Every pair of these forms consists of diametrally opposite points on the equator of corresponding tetrahedron edge.
There no exists another complex structure on the edges of tetrahedron.
Edge-spheres, which join this equatorial points of skew edges consist of integrable structures only.

''Maximum non-integrable'' structures on Lie group $SU(2)\times SU(2)$ take first sphere tangent vectors and replace its into the
second sphere tangent vector [7].
For example, structures $I_0$, $I_1$, $I_2$, $I_3$ are ''maximum non-integrable''.
The group $SO(3)\times SO(3)$ acts transitively with isotropy subgroup $diag(SO(3))$, the orbit of $I_0$
with respect to group $SO(3)\times SO(3)$ action is $SO(3)=SO(3)\times SO(3)/diag(SO(3))$.
Use formulas of corollary 2
one can find conditions, such that $J_{12}=J_{13}=J_{23}=J_{45}=J_{46}=J_{56}=0$. So we get solution
$$
\left\{
\begin{array}{l}
\cos(\varphi_1)=0,\\
\sin(\varphi_2)=0,\\
\sin(\varphi_3)=0
\end{array}
\right. \mbox{  or  }
\left\{
\begin{array}{l}
\sin(\varphi_1)=0,\\
\cos(\varphi_2)=0,\\
\cos(\varphi_3)=0
\end{array}
\right.
$$
The both solutions have the same form, so it is enough to consider only one conditions system.
Let
$$
M_+=\{\sin\psi(e^{14}+e^{25})+\cos\psi(e^{15}+e^{42})+e^{36}:-\pi\leq\psi\leq\pi\},
$$
$$
M_-=\{\sin\theta(e^{14}-e^{25})+\cos\theta(e^{15}-e^{42})+e^{36}:-\pi\leq\theta\leq\pi\}.
$$
are ''maximum non-integrable'' structures on the meridians of ${\mathcal E}_{03}$ and ${\mathcal E}_{12}$
accordingly.
Denote by $M_{p_+p_-}$ the half of big circle of meridian of general edge-sphere ${\mathcal E}_{p_-p_+}$, then
the set of all ''maximum non-integrable'' structures forms 3-dimensional subspace in $\mathbb{C}P^3$:
$$
\bigcup_
{\begin{array}{c}
p_-\in M_-\\
p_+\in M_+
\end{array}}M_{p_- p_+}
$$
 \vspace{1.5mm}
\textbf{Bibliography}
 \vspace{1.5mm}

[1] Buchshtaber V.M., Panov T.E. {\it Torus actions and their applications in topology and combinatorics. Univ.Lecture Series}, Vol. 24,
Amer. Math. Soc., Providence, RI, 2002, 152 p.

[2] Daurtseva N.A., Smolentsev N.K. {\it On the space of almost complex
structures}, Preprint -- arXiv: math.DG/0202139, (2002)

[3] Kobayashi S.,Nomizu K.:{\it Foundations of differential
geometry}, Vol.2, Intersciense Publishers, New York, London, (1969).

[4] Abbena, E. {\it Almost Hermitian geometry on six dimensional nilmanifolds}/E.~Abbena, S.~Garbiero, S.~Salamon// Ann. Scuola Norm. Sup. Pisa Cl.
Sci. {\bf 30} (2001), 147-–170.

[5] Calabi, E. {\it A class of compact complex manifolds which are not algebraic.}/E.~Calabi, B.~Eckmann//
Ann.Math. {\bf 58}(1935), 494--500.

[6] Daurtseva, N.A. {\it Invariant complex structures on $S^3\times S^3$}/N.A. Daurtseva//Electronic Journ. "Investigated
in Russia",2004, 888-893; (http://zhurnal.ape.relarn.ru/articles/2004/081e.pdf)

[7] Daurtseva, N.A. {\it Left-invariant almost nearly K$\mathrm{\ddot{a}}$hler structures on $SU(2)\times SU(2)$ in the tetrahedron visualization
for $\mathbb{C}P^3$}/N.~A.~Daurtseva//\\
arXiv:0608704[math.DG](2006) 12 p.

[8] Ivashkovich, S. {\it Complex curves in almost-complex manifolds and meromorphic hulls}/
S. Ivashkovich, V. Shevchishin//Bochum: Ruhr-Univ. Bochum, 1999. -- VI, 186 p. -- (Bochum Ruhr-Universitat: Schriftenreihe; h. 36)

[9] Magnin, L. {\it Left invariant complex structures on $U(2)$ and $SU(2)\times SU(2)$ revisited}/L. Magnin// preprint,
arXiv: 0809.1182 [math.RA] (2008) 25 p.

\end{document}